\documentclass[12pt]{article}
\usepackage{amsfonts}
\usepackage{a4wide}
\usepackage{amsmath,amscd,amsthm,a4,amssymb}

\begin{document}

\begin{center}
{\large \textbf{ON THE APPROXIMATION BY WEIGHTED\\[0pt]
RIDGE FUNCTIONS }}\footnote{%
Supported by INTAS Grant YSF-06-1000015-6283}

\bigskip

\bigskip

\textbf{Vugar E. Ismailov}

\smallskip

{\small Mathematics and Mechanics Institute}

{\small Azerbaijan National Academy of Sciences}

{\small Az-1141, Baku, Azerbaijan, e-mail: vugaris@mail.ru}
\end{center}

{\small \bigskip }

{\small \textbf{Abstract.} {We characterize the best $L_{2}$ approximation
to a multivariate function by linear combinations of ridge functions
multiplied by some fixed weight functions. In the special case when the
weight functions are constants, we propose explicit formulas for both the
best approximation and approximation error.} }


{\small \bigskip }

\begin{center}
{\small {\large \textbf{1. Introduction} } }
\end{center}

{\small A function $g\left( \mathbf{a}\cdot \mathbf{x}\right) ,$ where $%
\mathbf{a\in }\mathbb{R}^{n}\backslash \left\{ \mathbf{0}\right\} $, $%
\mathbf{x}\in \mathbb{R}^{n}$, $\mathbf{a}\cdot \mathbf{x}$ is the inner
product and $g$ is a univariate function, is called a \textit{ridge function
}(in $\mathbf{x}$) with the direction $\mathbf{a}$. These functions and
their linear combinations appear naturally in computerized tomography,
statistics, partial differential equations (where they are called \textit{%
plane waves}), neural networks, and approximation theory. Ridge
approximation in $L_{2}$ was actively studied in the late 90's by K.I.
Oskolkov [7], V.E. Maiorov [6], A. Pinkus [9], V.N. Temlyakov [10], P.
Petrushev [8] and others. }

{\small Let $D$ be the unit disk in $\mathbb{R}^{2}$. In [5], Logan and
Shepp along with other results gave a closed-form expression for the best $%
L_{2}$ approximation to a function $f\left( x_{1},x_{2}\right) \in
L_{2}\left( D\right) $ from the set }

{\small
\begin{equation*}
\mathcal{R}\left( \mathbf{a}^{1},...,\mathbf{a}^{r}\right) =\left\{
\sum\limits_{i=1}^{r}g_{i}\left( \mathbf{a}^{i}\cdot \mathbf{x}\right)
:g_{i}:\mathbb{R}\rightarrow \mathbb{R},~i=1,...,r\right\} .
\end{equation*}%
Their solution requires that the directions $\mathbf{a}^{1},...,\mathbf{a}%
^{r}$ be equally-spaced and involves finite sums of convolutions with
explicit kernels. In $n$ dimensional case, the author [3] obtained an
expression of simpler form for the best $L_{2}$ approximation to
square-integrable multivariate functions over some domain, provided that $%
r=n $ and the directions $\mathbf{a}^{1},...,\mathbf{a}^{r}$ are linearly
independent.}

{\small It should be noted that problems of approximation from the set $%
\mathcal{R}\left( \mathbf{a}^{1},...,\mathbf{a}^{r}\right) $ were also
considered in the uniform norm. For example, one essential approximation
method, its defects and advantages were discussed in [9]. Lin and Pinkus [4]
characterized $\mathcal{R}\left( \mathbf{a}^{1},...,\mathbf{a}^{r}\right) $,
i.e. they found means of determining if a continuous function $f$ (defined
on $\mathbb{R}^{n}$) is of the form $\sum\limits_{i=1}^{r}g_{i}\left(
\mathbf{a}^{i}\cdot \mathbf{x}\right) $ for some given $\mathbf{a}^{1},...,%
\mathbf{a^{r}}\in \mathbb{R}^{n}\backslash \left\{ \mathbf{0}\right\} $, but
unknown continuous $g_{1},...,g_{r}$. Two other characterizations of $%
\mathcal{R}\left( \mathbf{a}^{1},...,\mathbf{a}^{r}\right) $ may be found in
Diaconis and Shahshahani [2]. Buhmann and Pinkus [1] solved the inverse
problem: assume that we are given a function $f$ $\in $ $\mathcal{R}\left(
\mathbf{a}^{1},...,\mathbf{a}^{r}\right)$. How can we identify the functions
$g_{i},$ $i$ $=1,...,r$? }

{\small In this paper, we would like to consider the approximation from the
more general set }

{\small
\begin{equation*}
\mathcal{R}\left( \mathbf{a}^{1},...,\mathbf{a}^{r};~w_{1},...,w_{r}\right)
=\left\{ \sum\limits_{i=1}^{r}w_{i}(\mathbf{x})g_{i}\left( \mathbf{a}%
^{i}\cdot \mathbf{x}\right) :g_{i}:\mathbb{R}\rightarrow \mathbb{R}%
,~i=1,...,r\right\} ,
\end{equation*}%
where $w_{1},...,w_{r}$ are fixed multivariate functions. We are going to
characterize the best $L_{2}$ approximation in this set (see theorem 2.4)
for the case $r\leq n.$ Then, in the special case when the weight functions $%
w_{1},...,w_{r}$ are constants, we will prove two theorems on explicit
formulas for the best approximation and the error of approximation
respectively. Unfortunately, we do not yet know any reasonable answer to
these problems in other possible cases of $r.$ }

{\small \bigskip }

\begin{center}
{\small {\large \textbf{2. Characterization of the best approximation} } }
\end{center}

{\small Let $X$ be a subset of $\mathbb{R}^{n}$ with a finite Lebesgue
measure. Consider the approximation of a function $f\left( \mathbf{x}\right)
=f\left( x_{1},...,x_{n}\right) $ in $L_{2}\left( X\right) $ from the
manifold $\mathcal{R}\left( \mathbf{a}^{1},...,\mathbf{a}%
^{r};~w_{1},...,w_{r}\right) $, where $r\leq n.$ We suppose that the
functions $w_{i}(\mathbf{x})$ and the products $w_{i}(\mathbf{x})\cdot
g_{i}\left( \mathbf{a}^{i}\cdot \mathbf{x}\right) ,~i=1,...,r$, belong to
the space $L_{2}\left( X\right) .$ Besides, we assume that the vectors $%
\mathbf{a}^{1},...,\mathbf{a}^{r}$ are linearly independent. We say that a
function $g_{w}^{0}=\sum\limits_{i=1}^{r}w_{i}(\mathbf{x})g_{i}^{0}\left(
\mathbf{a}^{i}\cdot \mathbf{x}\right) $ in $\mathcal{R}\left( \mathbf{a}%
^{1},...,\mathbf{a}^{r};~w_{1},...,w_{r}\right) $ is the best approximation
(or extremal) to $f$ if }

{\small
\begin{equation*}
\left\Vert f-g_{w}^{0}\right\Vert _{L_{2}\left( X\right) }=\inf\limits_{g\in
\mathcal{R}\left( \mathbf{a}^{1},...,\mathbf{a}^{r};~w_{1},...,w_{r}\right)
}\left\Vert f-g\right\Vert _{L_{2}\left( X\right) }.
\end{equation*}
}

{\small Let the system of vectors $\{\mathbf{a}^{1},...,\mathbf{a}^{r},%
\mathbf{a}^{r+1},...,\mathbf{a}^{n}\}$ be a completion of the system $\{%
\mathbf{a}^{1},...,\mathbf{a}^{r}\}$ to a basis in $\mathbb{R}^{n}.$\ Let $%
J:X\rightarrow \mathbb{R}^{n}$ be the linear transformation given by the
formulas
\begin{equation*}
y_{i}=\mathbf{a}^{i}\cdot \mathbf{x,}\quad \,i=1,...,n.\eqno(2.1)
\end{equation*}%
Since the vectors $\mathbf{a}^{i},$ $i=1,...,n$, are linearly independent,
it is an injection. The Jacobian $\det J$ of this transformation is a
constant different from zero. }

{\small Let the formulas
\begin{equation*}
x_{i}=\mathbf{b}^{i}\cdot \mathbf{y},\;\;i=1,...,n,
\end{equation*}%
stand for the solution of linear equations (2.1) with respect to $%
x_{i},\;i=1,...,n.$ }

{\small Introduce the notation
\begin{equation*}
Y=J\left( X\right)
\end{equation*}%
and
\begin{equation*}
Y_{i}=\left\{ y_{i}\in \mathbb{R}:\;\;y_{i}=\mathbf{a}^{i}\cdot \mathbf{x}%
,\;\;\mathbf{x}\in X\right\} ,\,i=1,...,n.
\end{equation*}
}

{\small For any function $u\in L_{2}\left( X\right) ,$ put
\begin{equation*}
u^{\ast }=u^{\ast }\left( \mathbf{y}\right) \overset{def}{=}u\left( \mathbf{b%
}^{1}\cdot \mathbf{y},...,\mathbf{b}^{n}\cdot \mathbf{y}\right) .
\end{equation*}%
}

{\small It is obvious that $u^{\ast }\in L_{2}\left( Y\right) .$ Besides,
\begin{equation*}
\int\limits_{Y}u^{\ast }\left( \mathbf{y}\right) d\mathbf{y}=\left\vert \det
J\right\vert \cdot \int\limits_{X}u\left( \mathbf{x}\right) d\mathbf{x}\eqno%
(2.2)
\end{equation*}%
and
\begin{equation*}
\left\Vert u^{\ast }\right\Vert _{L_{2}\left( Y\right) }=\left\vert \det
J\right\vert ^{1/2}\cdot \left\Vert u\right\Vert _{L_{2}\left( X\right) }.%
\eqno(2.3)
\end{equation*}%
\qquad \qquad \qquad \qquad \qquad }

{\small Set }

{\small
\begin{equation*}
L_{2}^{i}=\{w_{i}^{\ast }(\mathbf{y})g\left( y_{i}\right) \in
L_{2}(Y)\},~i=1,...,r.
\end{equation*}
}

{\small We need the following auxiliary lemmas. }

{\small \textbf{Lemma 2.1.} \textit{Let $f\left( \mathbf{x}\right) \in
L_{2}\left( X\right) $. A function $\sum\limits_{i=1}^{r}w_{i}(\mathbf{x}%
)g_{i}^{0}\left( \mathbf{a}^{i}\cdot \mathbf{x}\right) $ is extremal to the
function $f\left( \mathbf{x}\right) $ if and only if \ $\sum%
\limits_{i=1}^{r}w_{i}^{\ast }(\mathbf{y})g_{i}^{0}\left( y_{i}\right) $ is
extremal from the space $L_{2}^{1}\mathit{\oplus }...\oplus L_{2}^{r}$ to
the function $f^{\ast }\left( \mathbf{y}\right) $.} }

{\small \bigskip }

{\small Due to (2.3) the proof of this lemma is obvious. }

{\small \textbf{Lemma 2.2.} \textit{Let $f\left( \mathbf{x}\right) \in
L_{2}\left( X\right) $. A function $\sum\limits_{i=1}^{r}w_{i}(\mathbf{x}%
)g_{i}^{0}\left( \mathbf{a}^{i}\cdot \mathbf{x}\right) $ is extremal to the
function $f\left( \mathbf{x}\right) $ if and only if}
\begin{equation*}
\int\limits_{X}\left( f\left( \mathbf{x}\right) -\sum\limits_{i=1}^{r}w_{i}(%
\mathbf{x})g_{i}^{0}\left( \mathbf{a}^{i}\cdot \mathbf{x}\right) \right)
w_{j}(\mathbf{x})h\left( \mathbf{a}^{j}\cdot \mathbf{x}\right) d\mathbf{x}%
=0\
\end{equation*}%
\textit{for any ridge function $h\left( \mathbf{a}^{j}\cdot \mathbf{x}%
\right) $ such that $\mathit{w}_{j}\mathit{(x)h}\left( \mathbf{a}^{j}\cdot
\mathbf{x}\right) $$\in L_{2}\left( X\right) \;\;j=1,...,r$.} }

{\small \bigskip }

{\small \textbf{Lemma 2.3.} \textit{The following formula is valid for the
error of approximation to a function $f\left( \mathbf{x}\right) $ in $%
L_{2}\left( X\right) $ from $\mathcal{R}\left( \mathbf{a}^{1},...,\mathbf{a}%
^{r};~w_{1},...,w_{r}\right) $:}
\begin{equation*}
E\left( f\right) =\left( \left\Vert f\left( \mathbf{x}\right) \right\Vert
_{L_{2}\left( X\right) }^{2}-\left\Vert \sum\limits_{i=1}^{r}w_{i}(\mathbf{x}%
)g_{i}^{0}\left( \mathbf{a}^{i}\cdot \mathbf{x}\right) \right\Vert
_{L_{2}\left( X\right) }^{2}\right) ^{\frac{1}{2}},
\end{equation*}%
\textit{where $\sum\limits_{i=1}^{r}w_{i}(\mathbf{x})g_{i}^{0}\left( \mathbf{%
a}^{i}\cdot \mathbf{x}\right) $ is the best approximation to $f\left(
\mathbf{x}\right) $.} }

{\small \bigskip }

{\small \ Lemmas 2.2 and 2.3 follow from the well-known facts of functional
analysis that the best approximation of an element $x$ in a Hilbert space $H$
from a linear subspace $Z$ of $H$ must be the image of $x$ via the
orthogonal projection onto $Z$ and the sum of squares of norms of orthogonal
vectors is equal to the square of the norm of their sum. }

{\small \bigskip }

{\small We say that $Y$ is an $r$-set if it can be represented as $%
Y_{1}\times ...\times Y_{r}\times Y_{0},$ where $Y_{0}$ is some set from the
space $\mathbb{R}^{n-r}.$ In special case, $Y_{0}$ may be equal to $%
Y_{r+1}\times ...\times Y_{n},$ but it is not necessary. By $Y^{\left(
i\right) },$ we denote the Cartesian product of the sets $%
Y_{1},...,Y_{r},Y_{0}$ except for $Y_{i},\;i=1,...,r$. That is, $Y^{\left(
i\right) }=Y_{1}\times ...\times Y_{i-1}\times Y_{i+1}\times ...\times
Y_{r}\times Y_{0},\,\ i=1,...,r$. }

{\small \bigskip }

{\small \textbf{Theorem 2.4.} \textit{Let $Y$ be an }$r$\textit{-set. A
function $\sum\limits_{i=1}^{r}w_{i}(\mathbf{x})g_{i}^{0}\left( \mathbf{a}%
^{i}\cdot \mathbf{x}\right) $ is the best approximation to $f(\mathbf{x)}$
if and only if}
\begin{equation*}
g_{j}^{0}\left( y_{j}\right) =\frac{1}{\int\limits_{Y^{\left( j\right)
}}w_{j}^{\ast 2}(\mathbf{y})d\mathbf{y}^{\left( j\right) }}%
\int\limits_{Y^{\left( j\right) }}\left( f^{\ast }\left( \mathbf{y}\right)
-\sum\limits_{\substack{ i=1  \\ i\neq j}}^{r}w_{i}^{\ast }(\mathbf{y}%
)g_{i}^{0}\left( y_{i}\right) \right) w_{j}^{\ast }(\mathbf{y})d\mathbf{y}%
^{\left( j\right) },\;\;j=1,...,r.\eqno(2.4)
\end{equation*}
}

{\small
\begin{proof} Necessity. Let a function $\sum\limits_{i=1}^{r}w_{i}(\mathbf{x}%
)g_{i}^{0}\left( \mathbf{a}^{i}\cdot \mathbf{x}\right) $ be extremal to $f$.
Then by lemma 2.1, the function $\sum\limits_{i=1}^{r}w_{i}^{\ast }(\mathbf{y%
})g_{i}^{0}\left( y_{i}\right) $ in $L_{2}^{1}\oplus ...\oplus L_{2}^{r}$ is
extremal to $f^{\ast }$. By lemma 2.2 and equality (2.2),
\begin{equation*}
\int\limits_{Y}f^{\ast }\left( \mathbf{y}\right) w_{j}^{\ast }(\mathbf{y}%
)h\left( y_{j}\right) d\mathbf{y}=\int\limits_{Y}w_{j}^{\ast }(\mathbf{y}%
)h\left( y_{j}\right) \sum\limits_{i=1}^{r}w_{i}^{\ast }(\mathbf{y}%
)g_{i}^{0}\left( y_{i}\right) d\mathbf{y}\eqno(2.5)
\end{equation*}%
for any product $w_{j}^{\ast }(\mathbf{y})h\left( y_{j}\right) $ in $%
L_{2}^{j},\;\;j=1,...,r$. Applying Fubini's theorem to the integrals in
(2.5), we obtain that
\begin{equation*}
\int\limits_{Y_{j}}h\left( y_{j}\right) \left[ \int\limits_{Y^{\left(
j\right) }}f^{\ast }\left( \mathbf{y}\right) w_{j}^{\ast }(\mathbf{y})d%
\mathbf{y}^{\left( j\right) }\right] dy_{j}=\int\limits_{Y_{j}}h\left(
y_{j}\right) \left[ \int\limits_{Y^{\left( j\right) }}w_{j}^{\ast }(\mathbf{y%
})\sum\limits_{i=1}^{r}w_{i}^{\ast }(\mathbf{y})g_{i}^{0}\left( y_{i}\right)
d\mathbf{y}^{\left( j\right) }\right] dy_{j}.
\end{equation*}%
Since $h\left( y_{j}\right) $ is an arbitrary function such that $%
w_{j}^{\ast }(\mathbf{y})h\left( y_{j}\right) \in L_{2}^{j}$,
\begin{equation*}
\int\limits_{Y^{\left( j\right) }}f^{\ast }\left( \mathbf{y}\right)
w_{j}^{\ast }(\mathbf{y})d\mathbf{y}^{(j)}=\int\limits_{Y^{\left( j\right)
}}w_{j}^{\ast }(\mathbf{y})\sum\limits_{i=1}^{r}w_{i}^{\ast }(\mathbf{y}%
)g_{i}^{0}\left( y_{i}\right) d\mathbf{y}^{\left( j\right) },\;\;j=1,...,r.
\end{equation*}%
Therefore,
\begin{equation*}
\int\limits_{Y^{\left( j\right) }}w_{j}^{\ast 2}(\mathbf{y})g_{j}^{0}\left( {%
y_{j}}\right) d\mathbf{y}^{\left( j\right) }=\int\limits_{Y^{\left( j\right)
}}\left( f^{\ast }\left( \mathbf{y}\right) -\sum\limits_{\substack{ i=1 \\ %
i\neq j}}^{r}w_{i}^{\ast }(\mathbf{y})g_{i}^{0}\left( y_{i}\right) \right)
w_{j}^{\ast }(\mathbf{y})d\mathbf{y}^{\left( j\right) },\;\;j=1,...,r.
\end{equation*}

Now, since $y_{j} \notin Y^{\left( j\right) }$, we obtain (2.4).

Sufficiency. Note that all the equalities in the proof of the necessity can
be obtained in the reverse order. Thus, (2.5) can be obtained from (2.4).
Then by (2.2) and lemma 2.2, we finally conclude that the function $%
\sum\limits_{i=1}^{r}w_{i}(\mathbf{x})g_{i}^{0}\left( \mathbf{a}^{i}\cdot
\mathbf{x}\right) $ is extremal to $f\left( \mathbf{x}\right) $.
\end{proof}
}

{\small In the following, $\left\vert Q\right\vert $ will denote the
Lebesgue measure of a measurable set $Q.$ The following corollary is
obvious. }

{\small \bigskip }

{\small \textbf{Corollary 2.5.} \textit{Let $Y$ be an }$r$\textit{-set. A
function $\sum\limits_{i=1}^{r}g_{i}^{0}\left( \mathbf{a}^{i}\cdot \mathbf{x}%
\right) $ in $\mathcal{R}\left( \mathbf{a}^{1},...,\mathbf{a}^{r}\right) $
is the best approximation to $f(\mathbf{x)}$ if and only if}
\begin{equation*}
g_{j}^{0}\left( y_{j}\right) =\frac{1}{\left\vert Y^{\left( j\right)
}\right\vert }\int\limits_{Y^{\left( j\right) }}\left( f^{\ast }\left(
\mathbf{y}\right) -\sum\limits_{\substack{ i=1  \\ i\neq j}}%
^{r}g_{i}^{0}\left( y_{i}\right) \right) d\mathbf{y}^{\left( j\right)
},\;\;j=1,...,r.
\end{equation*}
}

{\small \bigskip In [3], this corollary was proven for the case $r=n.$ }

{\small \bigskip }

\begin{center}
{\small {\large \textbf{3. Calculation of the approximation error}} }
\end{center}

{\small In this Section, we are going to establish explicit formulas for
both the best approximation and approximation error, provided that the
weight functions are constants. In this case, since we vary over $g_{i},$
the set $\mathcal{R}\left( \mathbf{a}^{1},...,\mathbf{a}%
^{r};~w_{1},...,w_{r}\right) $ coincide with $\mathcal{R}\left( \mathbf{a}%
^{1},...,\mathbf{a}^{r}\right) .$ Thus, without loss of generality, we may
assume that $w_{i}(\mathbf{x})=1$ for $i=1,...,r.$ }

{\small For brevity of the further exposition, introduce the notation }

{\small
\begin{equation*}
A=\int\limits_{Y}f^{\ast }\left( \mathbf{y}\right) d\mathbf{y}\text{ and \ }%
f_{i}^{\ast }=f_{i}^{\ast }(y_{i})=\int\limits_{Y^{\left( i\right) }}f^{\ast
}\left( \mathbf{y}\right) d\mathbf{y}^{\left( i\right) },~i=1,...,r.
\end{equation*}
}

{\small The following theorem is a generalization of the main result of [3]
from the case $r=n$ to the cases $r<n.$ }

{\small \textbf{Theorem 3.1.} \textit{Let $Y$ be an }$r$\textit{-set. Set
the functions}
\begin{equation*}
g_{1}^{0}\left( y_{1}\right) =\frac{1}{\left\vert Y^{\left( 1\right)
}\right\vert }f_{1}^{\ast }-\left( r-1\right) \frac{A}{\left\vert
Y\right\vert }
\end{equation*}%
\textit{and}
\begin{equation*}
g_{j}^{0}\left( y_{j}\right) =\frac{1}{\left\vert Y^{\left( j\right)
}\right\vert }f_{j}^{\ast },\;j=2,...,r.
\end{equation*}%
\textit{Then the function $\sum\limits_{i=1}^{r}g_{i}^{0}\left( \mathbf{a}%
^{i}\cdot \mathbf{x}\right) $ is the best approximation from $\mathcal{R}%
\left( \mathbf{a}^{1},...,\mathbf{a}^{r}\right) $ to $f\left( \mathbf{x}%
\right) $.} }

{\small The proof is the same as in [3]. It is sufficient to verify that the
functions $g_{j}^{0}\left( y_{j}\right) ,\;j=1,...,r$, satisfy the
conditions of corollary 2.5. This becomes obvious if note that
\begin{equation*}
\sum\limits_{\underset{i\neq j}{i=1}}^{r}\frac{1}{\left\vert Y^{\left(
j\right) }\right\vert }\frac{1}{\left\vert Y^{\left( i\right) }\right\vert }%
\int\limits_{Y^{\left( j\right) }}\left[ \int\limits_{Y^{\left( i\right)
}}f^{\ast }\left( \mathbf{y}\right) d\mathbf{y}^{\left( i\right) }\right] d%
\mathbf{y}^{\left( j\right) }=\left( r-1\right) \frac{1}{\left\vert
Y\right\vert }\int\limits_{Y}f^{\ast }\left( \mathbf{y}\right) d\mathbf{y}
\end{equation*}%
for $j=1,...,r$. }

{\small \bigskip }

{\small \textbf{Theorem 3.2.} \textit{Let $Y$ be an }$r$\textit{-set. Then
the error of approximation to a function $f(x)$ from the set $\mathcal{R}%
\left( \mathbf{a}^{1},...,\mathbf{a}^{r}\right) $ can be calculated by the
formula} }

{\small
\begin{equation*}
E(f)=\left\vert \det J\right\vert ^{-1/2}\left( \left\Vert f^{\ast
}\right\Vert _{L_{2}(Y)}^{2}-\sum_{i=1}^{r}\frac{1}{\left\vert Y^{\left(
i\right) }\right\vert ^{2}}\left\Vert f_{i}^{\ast }\right\Vert
_{L_{2}(Y)}^{2}+(r-1)\frac{A^{2}}{\left\vert Y\right\vert }\right) ^{1/2}.
\end{equation*}
}

{\small
\begin{proof} From Eq. (2.3), lemma 2.3 and theorem 3.1, it follows that

\begin{equation*}
E(f)=\left\vert \det J\right\vert ^{-1/2}\left( \left\Vert f^{\ast
}\right\Vert _{L_{2}(Y)}^{2}-I\right) ^{1/2},\eqno(3.1)
\end{equation*}%
where

\begin{equation*}
I=\left\Vert \sum_{i=1}^{r}\frac{1}{\left\vert Y^{\left( i\right)
}\right\vert }f_{i}^{\ast }-(r-1)\frac{A}{\left\vert Y\right\vert }%
\right\Vert _{L_{2}(Y)}^{2}.
\end{equation*}%
The integral $I$ can be written as a sum of the following four
integrals:

\begin{eqnarray*}
I_{1} &=&\sum_{i=1}^{r}\frac{1}{\left\vert Y^{\left( i\right) }\right\vert
^{2}}\left\Vert f_{i}^{\ast }\right\Vert
_{L_{2}(Y)}^{2},~I_{2}=\sum_{i=1}^{r}\sum\limits_{\substack{ j=1  \\ j\neq i
}}^{r}\frac{1}{\left\vert Y^{\left( i\right) }\right\vert }\frac{1}{%
\left\vert Y^{\left( j\right) }\right\vert }\int\limits_{Y}f_{i}^{\ast
}f_{j}^{\ast }d\mathbf{y,} \\
I_{3} &=&-2(r-1)\frac{1}{\left\vert Y\right\vert }A\sum_{i=1}^{r}\frac{1}{%
\left\vert Y^{\left( i\right) }\right\vert }\int\limits_{Y}f_{i}^{\ast }d%
\mathbf{y,}~I_{4}=(r-1)^{2}\frac{A^{2}}{\left\vert Y\right\vert }.
\end{eqnarray*}%
\qquad \qquad

It is not difficult to verify that

\begin{equation*}
\int\limits_{Y}f_{i}^{\ast }f_{j}^{\ast }d\mathbf{y=}\left\vert Y_{0}\times
\prod\limits_{\substack{ k=1  \\ k\neq i,j}}^{r}Y_{k}\right\vert A^{2},\text{
for }i,j=1,...,r,~i\neq j,\eqno(3.2)
\end{equation*}%
and

\begin{equation*}
\int\limits_{Y}f_{i}^{\ast }d\mathbf{y}=\left\vert Y_{0}\times \prod\limits
_{\substack{ k=1  \\ k\neq i}}^{r}Y_{k}\right\vert A,\text{ for }i=1,...,r.%
\eqno(3.3)
\end{equation*}%

Considering (3.2) and (3.3) in the expressions of $I_{2}$ and
$I_{3}$ respectively, we obtain that

\begin{equation*}
I_{2}=r(r-1)\frac{A^{2}}{\left\vert Y\right\vert }\text{ and }I_{3}=-2r(r-1)%
\frac{A^{2}}{\left\vert Y\right\vert }.
\end{equation*}%
Therefore,

\begin{equation*}
I=I_{1}+I_{2}+I_{3}+I_{4}=\sum_{i=1}^{r}\frac{1}{\left\vert Y^{\left(
i\right) }\right\vert ^{2}}\left\Vert f_{i}^{\ast }\right\Vert
_{L_{2}(Y)}^{2}-(r-1)\frac{A^{2}}{\left\vert Y\right\vert }.
\end{equation*}%
Now the last equality with (3.1) complete the proof.
\end{proof}
}

{\small \textbf{Example.} Consider the following set}

{\small
\begin{equation*}
X=\{\mathbf{x}\in \mathbb{R}^{4}:y_{i}=y_{i}(\mathbf{x})\in \lbrack
0;1],~i=1,...,4\},
\end{equation*}%
where}

{\small
\begin{equation*}
\left\{
\begin{array}{c}
y_{1}=x_{1}+x_{2}+x_{3}-x_{4} \\
y_{2}=x_{1}+x_{2}-x_{3}+x_{4} \\
y_{3}=x_{1}-x_{2}+x_{3}+x_{4} \\
y_{4}=-x_{1}+x_{2}+x_{3}+x_{4}%
\end{array}%
\right. \eqno(3.4)
\end{equation*}%
}

{\small \bigskip Let the function}

{\small
\begin{equation*}
f=8x_{1}x_{2}x_{3}x_{4}-\sum_{i=1}^{4}x_{i}^{4}+2\sum_{i=1}^{3}%
\sum_{j=i+1}^{4}x_{i}^{2}x_{j}^{2}
\end{equation*}%
be given over $X.$ Consider the approximation of this function from the set $%
\mathcal{R}\left( \mathbf{a}^{1},\mathbf{a}^{2},\mathbf{a}^{3}\right) ,%
\mathcal{\ }$where $\mathbf{a}^{1}=(1;1;1;-1),~\mathbf{a}^{2}=(1;1;-1;1),~%
\mathbf{a}^{3}=(1;-1;1;1).$ Putting $\mathbf{a}^{4}=(-1;1;1;1),$ we complete
the system of vectors $\mathbf{a}^{1},\mathbf{a}^{2},\mathbf{a}^{3}$ to the
basis $\{\mathbf{a}^{1},\mathbf{a}^{2},\mathbf{a}^{3},\mathbf{a}^{4}\}$ in $%
\mathbb{R}^{4}.$ The linear transformation $J$ defined by (3.4) maps the set
$X$ onto the set $Y=[0;1]^{4}.$ The inverse transformation is given by the
formulas}

{\small
\begin{equation*}
\left\{
\begin{array}{c}
x_{1}=\frac{1}{4}y_{1}+\frac{1}{4}y_{2}+\frac{1}{4}y_{3}-\frac{1}{4}y_{4} \\
x_{2}=\frac{1}{4}y_{1}+\frac{1}{4}y_{2}-\frac{1}{4}y_{3}+\frac{1}{4}y_{4} \\
x_{3}=\frac{1}{4}y_{1}-\frac{1}{4}y_{2}+\frac{1}{4}y_{3}+\frac{1}{4}y_{4} \\
x_{4}=-\frac{1}{4}y_{1}+\frac{1}{4}y_{2}+\frac{1}{4}y_{3}+\frac{1}{4}y_{4}%
\end{array}%
\right.
\end{equation*}%
}

{\small It can be easily verified that $f^{\ast }=y_{1}y_{2}y_{3}y_{4}$ and $%
Y$ is a $3$-set with $Y_{i}=[0;1],$ $i=1,2,3.$ Besides, $Y_{0}=[0;1].$ After
easy calculations we obtain that $A=\allowbreak \frac{1}{16};~$\ $%
f_{i}^{\ast }=\allowbreak \frac{1}{8}y_{i}$ for $i=1,2,3;$ $\det J=-16;$ $%
\left\Vert f^{\ast }\right\Vert _{L_{2}(Y)}^{2}=\frac{1}{81};$ $\left\Vert
f_{i}^{\ast }\right\Vert _{L_{2}(Y)}^{2}=\frac{1}{192},$ $i=1,2,3.$ Now from
theorems 3.1 and 3.2 it follows that the function $\frac{1}{8}%
\sum_{i=1}^{3}\left( \mathbf{a}^{i}\cdot \mathbf{x}\right) -\allowbreak
\frac{1}{8}$ is a best approximant from $\mathcal{R}\left( \mathbf{a}^{1},%
\mathbf{a}^{2},\mathbf{a}^{3}\right) $ to $f$ and $E(f)=\frac{1}{576}\sqrt{2}%
\sqrt{47}.$}

{\small \bigskip }

\end{document}